\documentclass{svmult}

\renewcommand{\email}[1]{\emailname: #1} 

\usepackage{mathptmx}       
\usepackage{helvet}         
\usepackage{courier}        

\usepackage{makeidx}         
\usepackage{graphicx}        
\usepackage[bottom]{footmisc}

\usepackage{latexsym}
\usepackage{amsmath}
\usepackage{amsfonts}
\usepackage{amssymb}
\usepackage{bm}

\usepackage{url}
\usepackage{algorithm}
\usepackage{algorithmic}
\usepackage[misc,geometry]{ifsym}

\spdefaulttheorem{assumption}{Assumption}{\upshape \bfseries}{\itshape}
\spdefaulttheorem{algo}{Algorithm}{\upshape \bfseries}{\itshape}

\newcommand{\bbN}{{\mathbb{N}}}

\newcommand{\bbR}{{\mathbb{R}}}
\newcommand{\bbS}{{\mathbb{S}}}

\DeclareSymbolFont{bbold}{U}{bbold}{m}{n}
\DeclareSymbolFontAlphabet{\mathbbold}{bbold}

\begin{document}

\title*{Tractability of Multivariate Problems
for Standard and Linear Information
in the Worst Case Setting: Part II}

\titlerunning{Standard Information in the Worst Case Setting}

\author{Erich Novak \and Henryk Wo\'zniakowski}
\authorrunning{E.~Novak, H.~Wo\'zniakowski}

\institute{
Erich Novak (\Letter)
\at Jena University, Math Institute, Ernst Abbe Platz 2, Jena, Germany \\
\email{erich.novak@uni-jena.de}
 \and
Henryk Wo\'zniakowski 
\at 
Department of Computer Science, Columbia University,
New York, NY 10027, USA 
\at 
Institute of Applied Mathematics, University of Warsaw, 
ul. Banacha 2, 02-097 Warszawa, Poland \\
\email{henryk@cs.columbia.edu} 
}


\maketitle

\index{Novak, Erich} 
\index{Wo\'zniakowski, Henryk} 

\paragraph{Dedicated to Ian H.~Sloan on the occassion of his 80th birthday.} 

\abstract{
We study QPT (quasi-polynomial tractability) 
in the worst case setting for
linear tensor product problems defined over Hilbert spaces.
We assume that the domain space 
is a reproducing kernel Hilbert space so that
function values are well defined.  
We prove QPT for algorithms that use only function values
under the three assumptions: 
\begin{enumerate}
\item
the minimal errors for the univariate case
decay polynomially fast to zero, 
\item 
the largest singular value for the univariate case is simple and 
\item
the eigenfunction 
corresponding to 
the largest singular value is a multiple of the function value at some point.
\end{enumerate} 
The first two assumptions are necessary for QPT. The third assumption
is necessary for QPT for some Hilbert spaces.
}

\section{Introduction}

In Part I \cite{NW14} we presented a lower error bound for approximating  linear 
multivariate operators defined over Hilbert spaces with algorithms that use 
function values. In this Part II we study upper bounds and algorithms 
for the same problem.  We want to understand 
the intrinsic difficulty of approximation of 
$d$-variate problems when $d$ is large. 
Algorithms that approximate $d$-variate problems may use
finitely many functionals from the class~$\Lambda^{\rm all}$
of information
or from the standard class~$\Lambda^{\rm std}$ of information.
The class~$\Lambda^{\rm all}$ consists of arbitrary linear functionals,
whereas the class~$\Lambda^{\rm std}$ consists of only function values.

We wish to approximate a $d$-variate problem in the worst case setting
to within an error threshold
$\varepsilon \in(0,1)$. The intrinsic difficulty is measured by the information 
complexity which is defined as the minimal number 
of linear functionals from the class 
$\Lambda\in\{\Lambda^{\rm all},\Lambda^{\rm std}\}$ 
which is needed 
to find an $\varepsilon $-approximation,
see \eqref{157157} for the precise definition.

Tractability deals with how the information complexity 
depends on $d$ and on $\varepsilon ^{-1}$, see
\cite{NW08,NW10,NW12}. In particular, we would like to know 
when the information complexity is exponential in $d$, the so-called
curse of dimensionality, and when 
we have a specific dependence on $d$ which is not exponential.
There are various ways of measuring the lack of exponential 
dependence and that leads to different notions of tractability.
In particular, we have polynomial tractability (PT) 
when the information complexity is polynomial in $d$ and $\varepsilon ^{-1}$,
and quasi-polynomial tractability (QPT) if the 
information complexity is at most proportional to 
$$
\exp\left(\,t\,(1+\ln\,\varepsilon ^{-1})(1+ \ln\,d)\right)
=\left(e\,\varepsilon^{-1}\right)^{t(1+\ln\,d)}
$$ 
for some non-negative $t$ independent of $d$ and $\varepsilon $.
This means that the exponent of~$\varepsilon^{-1}$ may depend weakly
on $d$ through $\ln\,d$.

In this paper we  study 
QPT for linear (unweighted) tensor product problems, $\bbS=\{S_d\}$
with $S_d=S_1^{\,\otimes\,d}$ and a compact linear non-zero $S_1:F_1\to G_1$ 
for Hilbert spaces $F_1$ and $G_1$. Since we want to use function values
we need to assume that $F_1$ is a reproducing kernel Hilbert space of univariate
functions defined on a non-empty~$D\subseteq\bbR$.
For simplicity we consider real valued functions.
By 
$$
K_1:D_1\times D_1\to\bbR
$$ 
we denote the reproducing kernel of $F_1$. 
Then $S^{\,\otimes\,d}:F_1^{\,\otimes\,d}\to G_1^{\,\otimes\,d}$ and
$F_1^{\,\otimes\,d}$ is a reproducing kernel Hilbert 
space of $d$-variate functions
defined on $D\times D\times\cdots\times D$ ($d$ times) 
with the reproducing kernel
$$
K_d(x,t)=\prod_{j=1}^dK_1(x_j,t_j)\ \ \ \ \
\mbox{for all}\ \ \ x_j,t_j\in D_1.
$$       

Obviously, tractability may depend 
on which class $\Lambda^{\rm std}$ or $\Lambda^{\rm all}$ is used.
Tractability results for $\Lambda^{\rm std}$ cannot be better than 
for $\Lambda^{\rm all}$.
The main question is when they are more or less the same. 
In particular, it is known when QPT holds for $\Lambda^{\rm all}$. Namely, let 
$\{\lambda_j,\eta_j\}$ be the ordered sequence of  
eigenvalues $\lambda_j$ and orthonormal 
eigenfunctions $\eta_j$ of $S_1^*S_1:F_1\to F_1$. 
Here $S_1^*: G_1 \to F_1$ is the adjoint operator of $S_1$. 
Let
$$
{\rm decay}_\lambda:=\sup\{\,r\ge0\, :\ \lim_{j\to\infty}
j^{\,r}\lambda_j=0\,\}
$$
denote the polynomial decay of the eigenvalues $\lambda_j$.
Since $S_1$ is assumed to be non-zero and compact, we have 
$\lambda_1>0$, $\lim_j\lambda_j=0$, and
${\rm decay}_\lambda$ is well defined.  However, it may happen that
${\rm decay}_\lambda=0$. 

It is known, see \cite{GW11}, that 
$$
\bbS \ \ \mbox{is \ QPT\  for \ $\Lambda^{\rm all}$}\ \ \  
\mbox{iff}\ \ \  
\lambda_2<\lambda_1\ \ \mbox{and}\ \ {\rm decay}_{\lambda}>0.
$$
 
Furthermore, if $\lambda_2>0$ 
then $\bbS$ is not PT for $\Lambda^{\rm all}$ (and for $\Lambda^{\rm std}$).
On the other hand, if $\lambda_2=\lambda_1>0$ then $\bbS$ suffers 
from the curse of dimensionality for the class $\Lambda^{\rm all}$ (and 
for $\Lambda^{\rm std}$).  

We now discuss QPT for $\Lambda^{\rm std}$. To motivate 
the need for the assumption
on the  eigenfunction $\eta_1$ corresponding 
to the largest eigenvalue $\lambda_1$,   
we cite a result from Part I, see
\cite{NW14}, for the Sobolev space 
$$
F_1\ \ \mbox{with the reproducing kernel}\ \ 
K_1^*(x,t)=1+\min(x,t)\ \ \mbox{for}\ \ x,t\in[0,1].
$$ 
\noindent
Then $\bbS$ suffers from the curse of dimensionality if
\begin{equation}\label{NWiffiff}
\eta_1\not=\pm\,[K_1^*(t,t)]^{-1/2}\,K_1^* (\cdot,t)=
\pm\,(1+t)^{-1/2}\,(1+\min(\cdot,t))\ \ \ \ 
\mbox{for all}\ \ \ t\in[0,1]. 
\end{equation}
Furthermore, for the approximation problem, 
$S_1f={\rm APP}_1f = f \in G_1=L_2([0,1])$, 
the assumption~\eqref{NWiffiff} holds, 
$\lambda_2<\lambda_1$ and ${\rm decay}_\lambda=2$.
Therefore for ${\rm APP} =\{ {\rm APP}_1^{\,\otimes\,d}\}$ we have 
$$
\mbox{Curse \ for\  $\Lambda^{\rm std}$ \ \ \ and \ \ \ QPT\  
\ for\  \ $\Lambda^{\rm all}$.}
$$ 
In this paper we prove that the assumption \eqref{NWiffiff} 
is essential for the curse  and QPT can hold for 
the class~$\Lambda^{\rm std}$ if \eqref{NWiffiff} is not satisfied.

This will be shown by establishing a result for general linear 
non-zero tensor product problems
for which $F_1$ is an arbitrary reproducing kernel Hilbert space with the 
reproducing kernel $K_1:D_1 \times D_1\to\bbR$. For the class
$\Lambda^{\rm std}$, the role of  
the sequence $\lambda=\{\lambda_j\}$ 
is replaced by the sequence $e=\{e_n(S_1)\}$ of the minimal worst case errors  
of algorithms that use at most $n$ function values. First of all, note that 
$$
\lim_ne_n(S_1)=0.
$$
Indeed, this holds for $S_1$ being a continuous linear functional,
see \cite{NW10} p. 79, and for a compact linear operator $S_1$ and for all positive $\varepsilon $
it is enough to approximate sufficiently well finitely many linear functionals.   
We define the polynomial decay of the minimal errors $e_n(S_1)$ as for the eigenvalues
by
$$
{\rm decay}_e:=\sup\{\,r\ge0\, :\ \lim_{n\to\infty}n^re_n(S_1)=0\,\}.
$$
The main result of this paper is the following theorem.

\begin{theorem}\label{NWmain}
\quad

Let $\bbS$ be a non-zero linear tensor product with a compact 
linear $S_1$ for which
\begin{itemize}
\item $\lambda_2<\lambda_1$,
\item ${\rm decay}_e>0$,
\item $\eta_1=\pm\,
K_1(t,t)^{-1/2}\,K_1(\cdot,t)$ \ for some $t\in D_1$.
\end{itemize}  
Then $\bbS$ is QPT for the class $\Lambda^{\rm std}$. 
\end{theorem}

We now comment on the assumptions of this theorem. The first
assumption is the same as for the class $\Lambda^{\rm all}$.
As already said, for $\lambda_2=\lambda_1>0$ we have 
the curse of dimensionality
for $\Lambda^{\rm std}$. 
The second assumption is necessary for QPT and the class $\Lambda^{\rm std}$. Indeed, if ${\rm decay}_e=0$ then even the univariate 
case cannot be solved polynomially in $\varepsilon ^{-1}$. This assumption corresponds to the assumption
${\rm decay}_\lambda>0$ for the class $\Lambda^{\rm all}$. For many problems we have ${\rm decay}_e={\rm decay}_\lambda$.
However, there are problems for which ${\rm decay}_\lambda=1$, ${\rm decay}_e=0$, and 
$e_n(S_1)$ can go to zero arbitrarily slowly, 
i.e., like $1/\ln(\ln(\cdots\ln(n))))$, where the number 
of $\ln$ can be arbitrarily large, see \cite{HNV} which is also reported 
in \cite{NW12} pp. 292-304. In this case, i.e., when ${\rm decay}_\lambda>0$ and ${\rm decay}_e=0$,
we have QPT for $\Lambda^{\rm all}$ and no QPT for $\Lambda^{\rm std}$.

We now discuss the last assumption which states that the eigenfunction 
$\eta_1$ corresponding to the largest
eigenvalue $\lambda_1$ is of a very special form. First of all, 
note that the scaling which is used 
in \eqref{NWiffiff} and here is needed 
to guarantee that $\Vert \eta_1 \Vert =1$. 
This implies that $K_1(t,t)>0$. 
For 
$\eta_1=\pm\,K_1(t,t)^{-1/2}\,K_1(\cdot,t)$  we have
$$
\left<  f,\eta_1\right>_{F_1}=\pm\,K_1(t,t)^{-1/2}\,
\left<  f, K_1(\cdot, t)\right>_{F_1}=
\pm\,K_1(t,t)^{-1/2}\,f(t).
$$
This means that the inner product $\left<  f,\eta_1\right>_{F_1}$ now can be computed exactly by one function value. 
Apparently, this important property allows us 
to achieve QPT for the class $\Lambda^{\rm std}$.
If this last assumption is not satisfied then we may decrease $F_1$ slightly by a rank 1 modification 
to obtain QPT for the modified problem, see Section 6. 

Theorem \ref{NWmain} will be proved constructively by presenting an algorithm $A_{d,\varepsilon }$ 
that computes an $\varepsilon $-approximation 
and uses at most $\mathcal{O}\left(\exp\left(t\,(1+\ln\,\varepsilon^{-1})(1+\ln\,d)\right)\right)$ function values for some $t$
independent of $d$ and $\varepsilon ^{-1}$.
The algorithm $A_{d,\varepsilon }$ is 
a modification of the Smolyak (sparse grid) algorithm applied to components 
of the operators $S_d$, see 
\cite{Smo63,WW99} and Chapter 15 of \cite{NW10} 
as well as Chapter 27 of \cite{NW12}.

It seems interesting to apply Theorem \ref{NWmain} to the space $F_1$ with the reproducing kernel $K_1^*$  
which was used before.
Combining the results of Part I with Theorem \ref{NWmain} we obtain the following corollary.

\begin{corollary}\label{NWsob}

 \quad

     Consider the spaces with $K_1^*$ as above. 
Then $\bbS$ is QPT for the class $\Lambda^{\rm std}$ iff  
\begin{itemize}
\item $\lambda_2<\lambda_1$,
\item ${\rm decay}_e>0$,
\item $\eta_1=\pm\,(1+t)^{-1/2}\,(1+\min(\cdot,t))$ \ for some $t\in D_1$.
\end{itemize}  
\end{corollary} 
 
\section{Preliminaries}

Let $S:F\to G$ be a continuous linear non-zero operator, 
where $F$ is a reproducing kernel Hilbert space
of real functions $f$ defined over 
a common non-empty domain~$D\subset \bbR^k$
for some positive integer~$k$, and $G$ is a Hilbert space.
We approximate $S$ by algorithms~$A_n$ that use at most~$n$ 
function values, i.e., we use the class $\Lambda^{\rm std}$. 
Without loss of generality we may assume that $A$ is linear,
see e.g., \cite{NW08,TWW88}. That is,
$$
A_nf=\sum_{j=1}^nf(t_j)\,g_j 
$$
for some $t_j\in D$ and $g_j\in S(F)\subseteq G$.
The worst case error of $A_n$ is defined as
$$
e(A_n)=\sup_{\|f\|_F\le1}\|Sf-A_nf\|_G=\|S-A_n\|_{F\to G}.
$$
For $n=0$, we take $A_n=0$ and then we obtain the initial error which is
$$
e(0)=e_0(S)=\|S\|_{F\to G}.
$$
Since $S$ is non-zero, the initial error is positive. 

We are ready to define the information complexity
for the class $\Lambda^{\rm std}$ and for the so-called normalized error criterion.  It is defined as
the minimal number of function values which are needed to reduce
the initial error by  a factor $\varepsilon \in(0,1)$. That is,  
\begin{equation}\label{157157}
n(\varepsilon ,S)=\min\{\,n\,:\ \exists\, A_n\ \mbox{such that}\ 
e(A_n)\le \varepsilon \,e_0(S)\}.
\end{equation}
Assume now that we have a sequence 
$$
\bbS=\{S_d\}_{d=1}^\infty
$$ 
of continuous linear non-zero operators 
$S_d:F_d\to G_d$, where $F_d$ is a reproducing kernel Hilbert space 
of real functions defined over a non-empty $D_d\subset \bbR^d$ and $G_d$ 
is a Hilbert space. 
In this case, we want to verify
how the information complexity $n(\varepsilon ,S_d)$ depends on $\varepsilon ^{-1}$
and $d$.  We say that 
$\bbS$ is quasi-polynomially tractable (QPT) for the class $\Lambda^{\rm std}$ 
iff there are non-negative numbers $C$ and $t$ such that
$$
n(\varepsilon ,S_d)\le C\,\exp\left(\,t\,
(1+\ln\,\varepsilon ^{-1})(1+\ln\,d)\right)\ \ \ 
\mbox{for all}\ \ \varepsilon \in(0,1), \ d \in \bbN .
$$
More about other tractability concepts 
can be found in \cite{NW08,NW10,NW12}.

\section{Linear Tensor Products}

We obtain a linear tensor product problem if the spaces $F=F_d$ and $G=G_d$
as well as $S=S_d$ are given by tensor products of $d$ copies of $F_1$
and $G_1$ as well as a continuous linear non-zero operator $S_1:F_1\to G_1$,
respectively, 
where $F_1$ is a reproducing kernel Hilbert space of real univariate 
functions defined over a non-empty
$D_1\subset \bbR$ and $G_1$ is a Hilbert space. 
To simplify the notation we assume that $F_1$ is of infinite dimension.
Then $F_d$ is an infinite dimensional space 
of~$d$-variate real functions defined 
on $D_d=D_1\times D_1\times\cdots\times D_1$ ($d$ times). 


We assume that $S_1$ is compact. Then all $S_d$ are also compact. 
Let $(\lambda_{j},\eta_{j})$ be the eigenpairs of $W_1=S_1^*S_1:F_1\to F_1$ 
with
$$
\lambda_1\ge\lambda_2\ge\dots\ge0\ \ \ \ \ \mbox{and}\ \ \ \ \ 
\left<  \eta_i,\eta_j\right>_{F_1}=\delta_{i,j}.
$$
Clearly, $\|S_1\|_{F_1\to G_1}=\sqrt{\lambda_1}$. 
Since $S$ is non-zero, $\lambda_1>0$. We have $f\in F_1$ iff 
$$
f=\sum_{j=1}^\infty\left<  f,\eta_j\right>_{F_1}\eta_j\ \ \ \mbox{with}\ \ \ 
\|f\|^2_{F_1}=\sum_{j=1}^\infty  \left<  f,\eta_j\right>_{F_1}^2<\infty. 
$$
Then
\begin{equation}\label{NWs1}
S_1f=\sum_{j=1}^\infty\left<  f,\eta_j\right>_{F_1}S_1\eta_j,
\end{equation}
where
$$
\left<  S_1\eta_i,S_1\eta_j\right>_{G_1}=\left< \eta_i,W_1\eta_j\right>_{F_1}=\lambda_j\delta_{i,j}.
$$
This means that the sequence $\{S_1\eta_j\}$ is orthogonal in $G_1$ and
$$
\|S_1f\|^2_{G_1}=\sum_{j=1}^\infty\left<  f,\eta_j\right>_{F_1}^2\lambda_j.
$$ 

For $d\ge2$, the eigenpairs $(\lambda_j,\eta_j)$ of $W_d=S_d^*S_d:F_d\to F_d$ 
are given in terms of the eigenpairs $(\lambda_j,\eta_j)$ of the univariate 
operator $W_1=S^*_1S_1:F_1\to F_1$. We have 
$$
\{\lambda_{d,j}\}_{j=1}^\infty=\{\lambda_{j_1}\lambda_{j_2}\cdots\lambda_{j_d}
\}_{j_1,j_2,\dots,j_d=1}^\infty.
$$
Similarly, the eigenfunctions of $W_d$ are of product form
$$
\{\eta_{d,j}\}_{j=1}^\infty=\{
\eta_{j_1}\otimes
\eta_{j_2}\otimes \cdots\otimes \eta_{j_d}
\}_{j_1,j_2,\dots,j_d=1}^\infty,   
$$
where
$$
[\eta_{j_1}\otimes
\eta_{j_2}\otimes \cdots\otimes \eta_{j_d}](x)=
\prod_{k=1}^d\eta_{j_k}(x_k)\ \ \ \ 
\mbox{for all}\ \ x=[x_1,\dots,x_d]\in D_d
$$
and
$$
\left<  \eta_{i_1}\otimes \eta_{i_2}\otimes\cdots\otimes \eta_{i_d},
\eta_{j_1}\otimes \eta_{j_2}\otimes\cdots\otimes \eta_{j_d}\right>_{F_d}=\delta_{i_1,j_1}
\delta_{i_2,j_2}\cdots\delta_{i_d,j_d}.
$$
Then $\|S_d\|_{F_d\to G_d}=\|W_d\|_{F_d\to F_d}^{1/2}=\lambda_1^{d/2}$.
Hence, the initial error is
$
e_0(S_d)=\lambda_1^{d/2}
$. 

We have $f\in F_d$ iff 
$$
f=
\sum_{(j_1,j_2,\dots,j_d)\in\bbN^d}\left<  f,\eta_{j_1}\otimes\cdots\otimes\eta_{j_d}\right>_{F_d}
\eta_{j_1}\otimes\cdots\otimes \eta_{j_d}
$$
with  
$$
\|f\|_{F_d}^2=
\sum_{(j_1,j_2,\dots,j_d)\in\bbN^d}\left<  f,\eta_{j_1}\otimes\cdots\otimes\eta_{j_d}\right>_{F_d}^2<\infty.
$$
In particular, for $x=(x_1,x_2,\dots,x_d)\in D_d$ we have
$$
f(x)=
\sum_{(j_1,j_2,\dots,j_d)\in\bbN^d}
\left<  f,\eta_{j_1}\otimes\cdots\otimes\eta_{j_d}\right>_{F_d}
\eta_{j_1}(x_1)\cdots\eta_{j_d}(x_d).
$$

\section{Decomposition of Linear Tensor Products}
In this section we assume, as in Theorem \ref{NWmain}, that
$$
\eta_1=\pm\,K_1(t,t)^{-1/2}K_1(\cdot,t) \ \ \ \mbox{for some\ \ $t\in D_1$}.
$$
Then for $j\ge2$ we obtain
$$
0=\left<  \eta_1,\eta_j\right>_{F_1}= K_1(t,t)^{-1/2}\,\eta_j(t).
$$
Hence, $\eta_j(t)=0$ for all $j\ge2$. This implies 
that 
$$
f(t,\dots,t)=\left<  f,\eta_1^{\,\otimes\,d}\right>_{F_d}\eta^d_1(t), 
$$
and for any $k=1,2,\dots, d-1$
and any vector $x=(t,\dots,t,x_{k+1},\dots,x_d)$ we have
\begin{equation}\label{NWx11}
f(x)=
\sum_{(j_{k+1},\dots,j_d)\in\bbN^{d-k}}
\left<  f,\eta_{1}^{\,\otimes\,k}\otimes\eta_{j_{k+1}}\cdots\otimes\eta_{j_d}\right>_{F_d}
[\eta_{1}(t)]^k\eta_{j_{k+1}}(x_{k+1})\cdots\eta_{j_d}(x_d).
\end{equation}

We start the 
decomposition of $S_d$ from the univariate case, $d=1$. {}From \eqref{NWs1} we have
$$
S_1=V_{1}+V_{2}
$$
with 
\begin{eqnarray*}
V_{1}f&=&\left<  f,\eta_1\right>_{F_1}S_1\eta_1= \pm K_1(t,t)^{-1/2}\,f(t)\,S_1\eta_1,\\
V_{2}f&=&\sum_{j=2}^\infty \left<  f,\eta_j\right>_{F_1}S_1\eta_j
\end{eqnarray*}
for all $f\in F_1$. Clearly, 
$$
\|V_{1}\|_{F_1\to G_1}=\|S_1\eta_1\|_{G_1}=\sqrt{\lambda_1}\ \ \ 
\mbox{and}\ \ \ \|V_{2}\|_{F_1\to G_1}=\|S_1\eta_2\|_{G_1}=\sqrt{\lambda_2}.
$$
We stress that we can compute $V_1f$ exactly by using one function value.

For $d\ge2$, we obtain
$$
S_d=(V_1+V_2)^{\,\otimes\,d}=\sum_{(j_1,j_2,\dots,j_d)\in\{1,2\}^d}
V_{j_1}\otimes V_{j_2}\otimes\cdots\otimes V_{j_d}.
$$

For $j=(j_1,j_2,\dots,j_d)\in\{1,2\}^d$ we define
$$
|j|_2=|\{j_i\,|\ j_i=2\}|
$$
as the number of indices equal to $2$. Clearly,
$$
\|V_{j_1}\otimes V_{j_2}\otimes\cdots\otimes V_{j_d}\|_{F_d\to G_d}
=\|V_1\|_{F_1\to G_1}^{\,d-|j|_2}\,\|V_2\|_{F_1\to G_1}^{\,|j|_2}=
\lambda_1^{\,(d-|j|_2)/2}\,\lambda_2^{\,|j|_2/2}.
$$

\section{Algorithms for Linear Tensor Products} 

We now derive an algorithm for linear tensor products for which
the assumptions of Theorem~\ref{NWmain} hold and we conclude
QPT for the class $\Lambda^{\rm std}$ from an estimate of the worst case error
of this algorithm.

To simplify the notation we assume that $\lambda_1=1$.
This can be done without loss of generality since otherwise we can replace
$S_1$ by $\lambda^{-1/2}_1S_1$. 

For $\lambda_1=1$ and due to the first assumption in Theorem \ref{NWmain},
we have 
$$
\|V_1\|_{F_1\to G_1}=1\ \ \ \mbox{and}\ \ \ 
\|V_2\|_{F_1\to G_1}=\lambda_2^{1/2}<1.
$$
Consider first $V_{2,d}=V_2^{\,\otimes d}$ with 
an exponentially small norm since 
$$
\|V_{2,d}\|_{F_d\to G_d}=\lambda_2^{d/2}.
$$
{}From the assumptions ${\rm decay}_e>0$ and $\lambda_2<1$,  
it was concluded in \cite{WW99}, see in particular
Lemma 1 and Theorem 2 of this paper, 
that for all $d\in\bbN$ there is a Smolyak/sparse grid algorithm
$$
A_{d,n}f=\sum_{m=1}^nf(t_{d,n,m})\,g_{d,n,m}\ \ \ 
\mbox{for all}\ \ \ f\in F_d
$$
for some $t_{d,n,m}\in D_d$ and $g_{d,n,m}=g_{d,n,m,1}\otimes\cdots\otimes
g_{d,n,m,d}$  with $g_{d,n,m,\ell}\in V_1(F_1)\subseteq G_1$, such that 
\begin{equation}\label{NWstrongpol}
e(A_{d,n})=\|V_{2,d}-A_{d,n}\|_{F_d\to G_d}\le \alpha\,n^{-r}\ \ \ 
\mbox{for all\ $d,n\in \bbN$}  
\end{equation}
for some positive $\alpha$ and $r$. We stress that 
$\alpha$ and $r$ are independent of $d$ and $n$. 

{}From the third assumption of Theorem~\ref{NWmain} we know that
$$
\eta_1=\delta\,K_1(t,t)^{-1/2}K_1(\cdot,t),\ \ \
\mbox{where}\ \ \ \delta\in\{-1,1\}.  
$$

For an integer $k\in[0,d]$, 
consider $V_1^{\,\otimes\,(d-k)}\,\otimes\,V_2^{\,\otimes\,k}$.
For $k=0$ we drop the second factor and for $k=d$ we drop the
first factor so that 
$V_1^{\,\otimes\,d}\,\otimes\,V_2^{\,\otimes\,0}= V_1^{\,\otimes\,d}$
and  $V_1^{\,\otimes\,0}\,\otimes\,V_2^{\,\otimes\,d}=V_2^{\,\otimes\,d}$. 

For $k=0$, we approximate $V_1^{\,\otimes\,d}$ by the algorithm
$$
A_{d,n,0}f=\frac{\delta^d}{K_1(t,t)^{d/2}}\,f(t,t,\dots,t)\,(S_1\eta_1)^{\,\otimes\,d}\ \ \ 
\mbox{for all} \ \ \ f\in F_d.
$$
Clearly, the error of this approximation is zero since
$A_{d,n,0}=V_1^{\,\otimes\,d}$ and $A_{d,n,0}$ 
uses one function value.  

For $k=d$, we approximate $V_2^{\,\otimes\,d}$ by the algorithm $A_{d,n}$ with error
at most $\alpha\,n^{-r}$.

For $k=1,2,\dots,d-1$, we approximate 
$V_1^{\,\otimes\,(d-k)}\otimes V_2^{\,\otimes k}$  by the algorithm
$$
A_{d,n,k}f=\frac{\delta^{d-k}}{[K_1(t,t)]^{(d-k)/2}}\,
\sum_{m=1}^nf(t,\dots,t,t_{k,n,m})
(S_1\eta_1)^{\,\otimes\,(d-k)}\,\otimes\,g_{k,n,m}
$$
for all $f\in F_d$. We now show that
\begin{equation}\label{NWlem1}
A_{d,n,k}=V_1^{\,\otimes\,(d-k)}\,\otimes\,A_{k,n}.
\end{equation}
Indeed, we know that $V_1\eta_j=0$ for all $j\ge2$. Then
\begin{eqnarray*}
&&(V_1^{\,\otimes \,(d-k)}\,\otimes\,A_{k,n})f\\
&=&
(V_1^{\,\otimes \,(d-k)}\,\otimes\,A_{k,n})\,\sum_{(j_1,\dots,j_d)\in\bbN^d}
\left<  f,\eta_{j_1}\otimes\cdots\otimes \eta_{j_d}\right>_{F_d}
\eta_{j_1}\otimes\cdots\otimes\eta_{j_d}\\
&=&\sum_{(j_1,\dots,j_d)\in\bbN^d}
\left<  f,\otimes_{\ell=1}^d\eta_{j_{\ell}}\right>_{F_d}
(V_1\eta_{j_1})\otimes \cdots\otimes (V_1\eta_{j_{d-k}}) 
\otimes A_{k,n}(\eta_{j_{d-k+1}}\otimes\cdots\otimes\eta_{j_d})
\\
&=&\alpha_k\,\sum_{(j_{d-k+1},\dots,j_d)\in\bbN^{k}}
\left<
  f,\eta_1^{\,\otimes\,(d-k)}\otimes_{\ell=1}^k\eta_{j_{d-k+\ell}}
\right>_{F_d}
\,(S_1\eta_1)^{\,\otimes\,(d-k)}\otimes\\
&&\qquad\quad\, 
\sum_{m=1}^n\left(\otimes_{\ell=1}^k\eta_{j_{d-k+\ell}}\right)(t_{k,n,m})
g_{k,n,m}\\
&=&\frac{\delta^{d-k}}{[K_1(t,t)]^{(d-k)/2}}
\,\sum_{m=1}^n h_m(S_1\eta_1)^{\,\otimes\,(d-k)}\,\otimes\,g_{k,n,m},
\end{eqnarray*}
where $\alpha_k=\delta^{d-k}K_1(t,t)^{(d-k)/2}\eta_1(t)^{d-k}$ and 
\begin{eqnarray*}
    h_m&=&\sum_{(j_{d-k+1},\dots,j_d)\in\bbN^k}
\left<  f,\eta_1^{\,\otimes\,(d-k)}\otimes\eta_{j_{d-k+1}} \otimes 
\cdots\otimes\eta_{j_d}\right>_{F_d}\\
&&\qquad\qquad\qquad\ \cdot\,\eta_1(t)^{d-k}
\left(\eta_{j_{d-k+1}}\otimes\cdots\otimes \eta_{j_d}\right)(t_{k,n,j}).
\end{eqnarray*}
{}From \eqref{NWx11} we conclude that
$$
h_j=f(t,\dots,t,t_{k,n,j})
$$
and 
$$
V_1^{\,\otimes \,(d-k)}\,\otimes\,A_{k,n}=A_{d,n,k},
$$
as claimed.  {}From this, we see that 
$$
V_1^{\,\otimes \,(d-k)}\,\otimes\,V_2^{\,\otimes\,k}-A_{d,n,k}=
V_1^{\,\otimes \,(d-k)}\,\otimes\,(V_2^{\,\otimes\,k}-A_{k,n})
$$
and
$$
e(A_{d,n,k})=\|V_2^{\,\otimes\,k}-A_{k,n}\|_{F_k\to G_k}\le \alpha\,n^{-r}.
$$
\vskip 1pc
We now explain how we approximate
$V_{j_1}\otimes\cdots\otimes V_{j_d}$ for an arbitrary 
$$
j=(j_1,\dots,j_d)\in \{1,2\}^d.
$$ 
The idea is the same as before, i.e.,
for the indices $j_\ell=1$ we approximate $V_{j_\ell}$ by 
itself, and for the rest of the indices, which are equal to $2$, 
we apply the Smolyak/sparse grid algorithm for proper parameters. 
More precisely, let $k=|j|_2$. The cases $k=0$ and $k=d$ have been 
already considered. Assume then that $k\in[1,d-1]:= \{ 1, 2, \dots , d-1 \}$. 
Let $\ell_i\in[1,d]$ be the $ith$ occurrence of $2$ in the vector $j$,
i.e., $1\le \ell_1 <\ell_2<\cdots<\ell_k\le d$, and 
$j_{\ell_1}=j_{\ell_2}=\dots=j_{\ell_k}=2$.  

Define the algorithm 
$$
A_{d,n,j}f=\frac{\delta^{d-k}}{K_1(t,t)^{(d-k)/2}}\,
\sum_{m=1}^nf(y_{d,n,j,m})\,h_{d,n,j,m,1}\otimes\cdots\otimes
h_{d,n,j,m,d},
$$
where the vector $y_{d,n,j,m}=(y_{d,n,j,m,1},\dots,y_{d,n,j,m,d})$ 
is given by
$$
y_{d,n,j,m,\ell}=
\begin{cases}
t&\ \ \mbox{if\ } j_\ell=1,\\
t_{k,n,m,i}&\ \ \mbox{if\ } j_\ell=2\ \mbox{and}\ \ell=\ell_i,
\end{cases} 
$$
and 
$$
h_{d,n,j,m.\ell}=
\begin{cases}
S_1\eta_1&\ \ \mbox{if\ } j_\ell=1,\\
g_{k,n,m,i}&\ \ \mbox{if\ } j_\ell=2 \ \mbox{and}\ \ell=\ell_i
\end{cases} 
$$
for $\ell=1,2,\dots,d$. 

The error of the algorithm $A_{d,n,j}$ is the same 
as the error of the algorithm $A_{d,n,|j|_2}$ 
since for (unweighted tensor) products 
the permutation of indices does not matter. 

Hence, for all $j\in\{1,2\}$, 
the algorithm $A_{d,n,j}$ uses at most $n$ function values and
\begin{equation}\label{NWerrorj}
e(A_{d,n,j})\le\alpha\,n^{-r},
\end{equation}
and this holds for all $d$. 
\vskip 1pc
We now define an algorithm which approximates $S_d$ with error 
at most $\varepsilon \in(0,1)$.
The idea of this algorithm is based on approximation of all
$V_{j_1}\otimes \cdots \otimes V_{j_d}$ 
whose norm is $\|V_2\|^{|j|_2}=\lambda_2^{|j|_2/2}$.
If $\lambda_2^{|j|_2/2}\le \varepsilon /2$ we approximate
$V_{j_1}\otimes \cdots \otimes V_{j_d}$ 
by zero otherwise by the algorithm $A_{d,n,j}$  for specially chosen $n$.
More precisely, let
$$
k=\min\left(d,\left\lceil 
\frac{2\ln\,\frac2{\varepsilon }}{\ln\,\frac1{\lambda_2}}\right\rceil\right).
$$ 
Define the algorithm
\begin{equation}\label{NWalgeps}
A_{d,n,\varepsilon }=\sum_{j\in \{1,2\}^d}A_{d,n,\varepsilon,j}
\end{equation}
with
$$
A_{d,n,\varepsilon,j}=\begin{cases}0&\ \ \ \mbox{if}\ |j|_2>k,\\
A_{d,n,j}&\ \ \ \mbox{if}\ |j|_2\le k.
\end{cases}
$$
Note that non-zero terms in \eqref{NWalgeps} correspond to 
$|j|_2\le k$ and each of them uses
at most $n$ function values. Therefore the algorithm $A_{d,n,\varepsilon }$ uses at most 
$$
{\rm card}(A_{d,n,\varepsilon })\le n\,\sum_{\ell=0}^{k}\binom{d}{\ell}
$$
function values.

We now analyze the error of $A_{d,n,\varepsilon }$. We have
$$
S_d-A_{d,n,\varepsilon }=\sum_{j\in \{1,2\}^d,\ |j|_2\le k}\left(
V_{j_1}\otimes\cdots\otimes V_{j_d}
-A_{d,n,\varepsilon,    j}\right) +
\sum_{j\in \{1,2\}^d,\ |j|_2> k}
V_{j_1}\otimes\cdots\otimes V_{j_d}.
$$
Note that the second operator in the sum above is zero if $k=d$.
For $k<d$ the terms of the second operator are orthogonal
and therefore it has norm at most $\lambda_2^{k/2}\le \varepsilon /2$
by the definition of $k$. 

{}From \eqref{NWerrorj} we conclude
$$
\|S_d-A_{d,n,\varepsilon }\|_{F_d\to G_d}\le \alpha\,n^{-r}\,\sum_{\ell=0}^{k}\binom{d}{\ell}
\ +\ \varepsilon /2.
$$
We now consider two cases $k\le d/2$ and $k>d/2$. 
We opt for simplicity at the expense of some error
overestimates which are still enough to establish QPT. 
\begin{itemize}
\item Case $k\le d/2$.

Then the binomial coefficients $\binom{d}{\ell}$ are increasing and 
$$
\sum_{\ell=0}^{k}\binom{d}{\ell}\le (k+1)\,
\binom{d}{k}\le(k+1)\,\frac{d^k}{k!}\le 2d^k.
$$
If we take $n$ such that
\begin{equation}\label{NWdefn}
\frac{2\alpha\,d^k}{n^r}\le\varepsilon /2
\end{equation}
then 
$$
e(A_{d,n,\varepsilon })=\|S_d-A_{d,n,\varepsilon }\|_{F_d\to G_d}\le \varepsilon .
$$
Since $k\le 1+2\ln(2\varepsilon ^{-1})/\ln(\lambda_2^{-1})$, we have 
$$
d^k\le \alpha_1(1+\varepsilon ^{-1})^{\alpha_2\,(1+\ln\,d)}
$$
for some $\alpha_1$ and $\alpha_2$ independent of $d$ and $\varepsilon ^{-1}$. 
Therefore
$$
n=\mathcal{O}\left(\exp\left(\mathcal{O}((1+\ln\,\varepsilon ^{-1})(1+\ln\,d))\right)\right)
$$
satisfies \eqref{NWdefn}. Furthermore, the cardinality of $A_{d,n,\varepsilon }$ is bounded 
by 
$$
2d^k\,n=
\mathcal{O}\left(\exp\left(\mathcal{O}((1+\ln\,\varepsilon ^{-1})(1+\ln\,d))\right)\right).
$$
\item Case $k>d/2$.

We now have $d\le 2k\le 2(1+2\ln(2\varepsilon ^{-1})/\ln(\lambda_2^{-1}))
=\mathcal{O}(1+\ln\,\varepsilon ^{-1})$.
We estimate $\sum_{\ell=0}^k\binom{d}{\ell}$ by $2^d=\exp(\mathcal{O}(1+\ln\,\varepsilon ^{-1}))$.
Then $2\alpha\,2^d\,n^{-r}\le \varepsilon /2$ for
$$
n=\mathcal{O}\left(\exp\left(\mathcal{O}(1+\ln\,\varepsilon ^{-1})\right)\right). 
$$
Hence
$$
e(A_{d,n,\varepsilon })\le \varepsilon 
$$
and the cardinality of $A_{d,n,\varepsilon }$ is bounded by
$$
2^d\,n=\mathcal{O}\left(\exp\left(\mathcal{O}(1+\ln\,\varepsilon ^{-1})\right)\right).
$$
\end{itemize}
\vskip 1pc
In both cases, $k\le d/2$ and $k>d/2$, we show that the error of the algorithm $A_{d,n,\varepsilon }$
is at most $\varepsilon $ and the number of function values used by this algorithm is at most
$$
\alpha_3\exp\left(\alpha_4\,(1+\ln\,\varepsilon ^{-1})(1+\ln\,d)\right)
$$
for some $\alpha_3$ and $\alpha_4$ independent of $\varepsilon ^{-1}$ and $d$. This shows that
the problem $\bbS=\{S_d\}$ is QPT. This also proves Theorem \ref{NWmain}. 

\section{Final Comments}

Let us assume, as in Theorem~1, that 
$\bbS$ is a non-zero linear tensor product problem 
with a compact linear $S_1$ for which
\begin{itemize}
\item $\lambda_2<\lambda_1$,
\item ${\rm decay}_e>0$, 
\end{itemize}  
but the last condition is not fulfilled, i.e., 
$$
\eta_1\not=\pm\,K_1(t,t)^{-1/2}\,K_1(\cdot,t)\ \ \ 
\mbox{for all}\ \ \ t\in D.
$$ 
Then, as we have seen, we cannot in general conclude
QPT for the class $\Lambda^{\rm std}$. 

We can ask whether we can modify the problem somehow, 
by decreasing the class $F_1$, in order to obtain QPT for the 
modified (smaller) spaces. It turns out that 
this is possible. 
For notational convenience 
we assume again that $\lambda_1 =1$. 

Since $\eta_1$ is non-zero, there exists a point $t^*\in D$ such that
$\eta_1(t^*)\not=0$. Define 
$$
\widetilde F_1 =  \{ f\in F_1 \mid \left< f,\eta_1\right>_{F_1} =
[\eta_1(t^*)]^{-1} f(t^*) \} .
$$
Note that $\eta_1\in\widetilde F_1$ and $\widetilde F_1$ is a linear
subspace of $F_1$. Let
$$
\widetilde f=\eta_1-\frac{K_1(\cdot,t^*)}{\eta_1(t^*)}.
$$
Clearly, $\widetilde f\in F_1$ and $\widetilde f\not=0$. 
Then $\widetilde F_1$ can be rewritten as 
$$
\widetilde F_1 =  \{ f\in F_1 \mid \big< f,\widetilde f\big>_{F_1} =0\}.
$$ 
It is easy to verify that the reproducing kernel $\widetilde K_1$
of $\widetilde F$ is
$$
\widetilde K_1(x,y)=K_1(x,y)-\frac{\widetilde f(x)\,\widetilde f(y)}
{\|\widetilde f\|^{\,2}}\ \ \ \ \ \mbox{for all}\ \ \ x,y\in D.
$$
Furthermore, it is also easy to check that
$$
\eta_1=\widetilde K_1(t^*,t^*)^{-1/2}\,\widetilde K_1(\cdot,t^*).
$$
The operator $\widetilde S_1=S_1\big|_{\widetilde F_1}$,
which is the restriction of $S_1$ to the subspace $\widetilde F_1$
satisfies  all assumptions of Theorem~1. 
Indeed, let $\widetilde \lambda_n$ be the ordered eigenvalues 
of 
$$
\widetilde W_1=\widetilde S_1^*\widetilde S_1:\widetilde F_1
\to \widetilde F_1.
$$
Since $\eta_1\in\widetilde F_1$ we have 
$\widetilde \lambda_1=\lambda_1=1$, whereas $\widetilde \lambda_n
\le \lambda_n$ for all $n\ge2$ since $\widetilde F_1\subseteq F_1$.
Therefore $\widetilde \lambda_2<\widetilde\lambda_1$.
Similarly, for     both classes $\Lambda^{\rm all}$ and $\Lambda^{\rm std}$,  
the minimal worst case errors for $\widetilde S_1$ are no larger than
the minimal worst case errors for $S_1$
Hence, applying Theorem 1 for $\widetilde S_1^{\otimes \,d}$,
we conclude QPT for the modified problem.

\end{document}